\documentclass[reqno,dvips]{amsart}

\usepackage{mathrsfs}
\usepackage{amsmath}
\usepackage{amssymb}
\usepackage{cite}
\usepackage{latexsym}
\usepackage{graphicx}

\usepackage{color}

\theoremstyle{plain}
\newtheorem{thm}{Theorem}[section]

\newtheorem{dfn}[thm]{Definition}
\newtheorem{prop}[thm]{Proposition}
\newtheorem{rmk}[thm]{Remark}

\def\D{\mathrm{D}}
\def\I{\mathscr{I}}

\def\R{\mathrm{R}}
\def\T{\mathscr{T}}

\def\d{\mathrm{d}}
\def\h{\mathrm{h}}

\def\Cset{\mathbb{C}}
\def\Nset{\mathbb{N}}
\def\Rset{\mathbb{R}}
\def\Sset{\mathbb{S}}
\def\Zset{\mathbb{Z}}

\def\epsilon{\varepsilon}

\DeclareMathOperator{\sech}{sech}

\DeclareMathOperator{\sn}{sn}
\DeclareMathOperator{\cn}{cn}
\DeclareMathOperator{\dn}{dn}

\makeatletter
 \@addtoreset{equation}{section}
\makeatother
\def\theequation{\arabic{section}.\arabic{equation}}

\begin{document}

% **********************************************************
% Title of This Paper
% **********************************************************

\title[Nonintegrability of Forced Nonlinear Oscillators]%
{Nonintegrability of Forced Nonlinear Oscillators}
\thanks{This work was partially supported by the JSPS KAKENHI Grant Numbers
 JP17H02859 and JP19J22791 .}

\author{Shoya Motonaga}
\author{Kazuyuki Yagasaki}

\address{Department of Applied Mathematics and Physics, Graduate School of Informatics,
Kyoto University, Yoshida-Honmachi, Sakyo-ku, Kyoto 606-8501, JAPAN}
\email{mnaga@amp.i.kyoto-u.ac.jp}
\email{yagasaki@amp.i.kyoto-u.ac.jp}

\date{\today}

\begin{abstract}
In recent papers by the authors
 (S.~Motonaga and K.~Yagasaki,
 Obstructions to integrability of nearly integrable dynamical systems near regular level sets, 
 submitted for publication, and
 K.~Yagasaki,
 Nonintegrability of nearly integrable dynamical systems near resonant periodic orbits,
 submitted for publication),
 two different techniques which allow us
 to prove the real-analytic or complex-meromorphic nonintegrability of forced nonlinear oscillators
 having the form of time-periodic perturbations of single-degree-of-freedom Hamiltonian systems
 were provided.
Here the concept of nonintegrability in the Bogoyavlenskij sense is adopted
 and the first integrals and commutative vector fields are also required to depend
 real-analytically or complex-meromorphically on the small parameter.
In this paper we review the theories and continue to demonstrate their usefulness.
In particular, we consider the periodically forced damped pendulum
 and prove its nonintegrability in the above meaning.
\end{abstract}
\maketitle

% **********************************************************
% Section 1
% **********************************************************

\section{Introduction}

In this paper we consider two-dimensional nonlinear systems of the form
\begin{equation}
\dot{x}=J\D H(x)+\epsilon g(x,\omega t),\quad
x\in\Rset^2,
\label{eqn:syse}
\end{equation}
where $\epsilon$ is a small parameter such that $0<|\epsilon|\ll 1$,
 $\omega>0$ is a constant, $H:\Rset^2\to\Rset$ and $g:\Rset^2\times\Sset^1$ are analytic,
 and $J$ is the $2\times 2$ symplectic matrix,
\[
J=\begin{pmatrix}
0 & 1\\
-1 & 0
\end{pmatrix}.
\]
When $\epsilon=0$, Eq.~\eqref{eqn:syse} becomes a planar Hamiltonian system
\begin{equation}
\dot{x}=J\D_x H(x)
\label{eqn:sys0}
\end{equation}
with a Hamiltonian function $H(x)$.

Systems of the form \eqref{eqn:syse} represent many forced nonlinear oscillators
 and have attracted much attention \cite{GH83,W90}.
In particular, perturbation techniques called the homoclinic and subharmonic Melnikov methods 
 have been developed:
The homocliinic Melnikov method
 enables us to discuss the existence of transverse homoclinic orbits and their bifurcations
 \cite{GH83,M63,W90},
 and the subhamonic Melnikov method 
 to discuss the existence of periodic orbits and their stability and bifurcations
 \cite{GH83a,GH83,W90,Y96,Y02,Y03}.
The techniques have been successfully applied
 to reveal the dynamics of numerous forced nonlinear oscillators.
See \cite{GH83a,GH83,W90,Y96,Y02,Y03} for the details.
Many of these nonlinear oscillators have been believed to be nonintegrable for $\epsilon\neq 0$
 although no mathematical proof has been given.
Moreover, they are planar Hamiltonian systems and consequently integrable
 in the Liouvilee sense \cite{A89,M99} when $\epsilon=0$.

In recent two papers \cite{MY21b,Y21b}
 the nonintegrability of \eqref{eqn:syse} in the following Bogoyavlenskij sense \cite{B98}
 were studied.
 
\begin{dfn}[Bogoyavlenskij]
\label{dfn:1a}
For $n,q\in\Nset$ such that $1\le q\le n$,
 an $n$-dimensional dynamical system
\[
\dot{x}=f(x),\quad x\in\Rset^n\text{ or }\Cset^n,
\]
is called \emph{$(q,n-q)$-integrable} or simply \emph{integrable} 
 if there exist $q$ vector fields $f_1(x)(:= f(x)),f_2(x),\dots,f_q(x)$
 and $n-q$ scalar-valued functions $F_1(x),\dots,F_{n-q}(x)$ such that
 the following two conditions hold:
\begin{enumerate}
\setlength{\leftskip}{-1.8em}
\item[\rm(i)]
$f_1(x),\dots,f_q(x)$ are linearly independent almost everywhere and commute with each other,
 i.e., $[f_j,f_k](x):=\D f_k(x)f_j(x)-\D f_j(x)f_k(x)\equiv 0$ for $j,k=1,\ldots,q$,
 where $[\cdot,\cdot]$ denotes the Lie bracket$;$
\item[\rm(ii)]
The derivatives $\D F_1(x),\dots, \D F_{n-q}(x)$ are linearly independent almost everywhere
 and $F_1(x),\dots,F_{n-q}(x)$ are first integrals of $f_1, \dots,f_q$,
 i.e., $\D F_k(x)\cdot f_j(x)\equiv 0$ for $j=1,\ldots,q$ and $k=1,\ldots,n-q$,
 where ``$\cdot$'' represents the inner product.
\end{enumerate}
We say that the system is \emph{analytically} $($resp. \emph{meromorphically}$)$ \emph{integrable}
 if the first integrals and commutative vector fields are analytic $($resp. meromorphic$)$. 
\end{dfn}
Definition~\ref{dfn:1a} is considered as a generalization of %complete
 Liouville-integrability for Hamiltonian systems \cite{A89,M99}
 since an $n$-degree-of-freedom Liouville-integrable Hamiltonian system with $n\ge 1$
 has not only $n$ functionally independent first integrals
 but also $n$ linearly independent commutative (Hamiltonian) vector fields
 generated by the first integrals.

In \cite{MY21b} a technique which allows us to prove the real-analytic nonintegrability
 of nearly integrable dynamical systems containing \eqref{eqn:syse} was developed,
  is based on the results of \cite{MY21a}.
It is also considered as an extension of the classical results
 of Poincar\'e \cite{P92} and Kozlov \cite{K83,K96} (see also \cite{MY21b,Y21b}).
On the other hand, in \cite{Y21a}
 a different technique for complex-meromorphic nonintegrability
 of nearly integrable dynamical systems based on versions due to Ayoul and Zung \cite{AZ10}
 of the Morales-Ramis and Morales-Ramis-Sim\'o theories \cite{M99,MR01,MRS07} 
 was developed and applied to \eqref{eqn:syse} in \cite{Y21b}.
It was successfully applied to the restricted three-body problem in \cite{Y21a,Y21c}.
In both the theories of \cite{MY21b,Y21b} for \eqref{eqn:syse},
 the first integrals and commutative vector fields are required to depend
 real-analytically or complex-meromorphically on the small parameter $\epsilon$
 near $\epsilon=0$.
Moreover, they were applied to the periodically forced damped Duffing oscillator \cite{D18,H79,U78}
 and their usefulness was demonstrated there.

Here we review the theories of \cite{MY21b,Y21b} for \eqref{eqn:syse}
 and continue to demonstrate their usefulness.
In particular, we consider the periodically forced damped pendulum \cite{GH83,H99}
\begin{equation}
\dot{x}_1=x_2,\quad
\dot{x}_2=-\sin x_1+\epsilon(\beta\cos\omega t-\delta x_2),\quad
(x_1,x_2)\in\Sset^1\times\Rset,
\label{eqn:fdp}
\end{equation}
and apply the two theories to it,
 where $\beta,\delta\ge 0$ are constants and $\Sset^1=\Rset/2\pi\Zset$.
The system \eqref{eqn:fdp} also provides a mathematical model
 for ac driven Josephson junctions \cite{KM86,HH86,SS85}.
Complicated dynamics was observed in numerical simulations of \eqref{eqn:fdp}.
See, e.g., \cite{H99}.
For $0<|\epsilon|\ll 1$, the existence of chaos in \eqref{eqn:fdp} when
\begin{equation}
\frac{\beta}{\delta}>\frac{4}{\pi}\cosh(\tfrac{1}{2}\pi\omega)
\label{eqn:chaos}
\end{equation}
was proved by the homoclinic Melnikov method \cite{GH83,M63,W90} in \cite{HH86,SS85}.
Moreover, a rigorous numerical proof of chaos for $\epsilon\beta=1$ and $\epsilon\delta=0.1$
 was given in \cite{BCGH08}.
Thus, the system \eqref{eqn:fdp} is believed to be nonintegrable,
 but its mathematical proof has not been given yet, to the authors' knowledge.
We prove the Bogoyavlenskij -nonintegrability of \eqref{eqn:fdp}
 when the first integrals and commutative vector fields are also required to depend
 real-analytically or complex-meromorphically on $\epsilon$ near $\epsilon=0$.

The outline of this paper is as follows:
In Sections~2 and 3, respectively,
 we briefly describe the theories of \cite{MY21b} and \cite{Y21b} for \eqref{eqn:syse}.
We apply the two theories to discuss the Bogoyavlenskij-nonintegrability
 of the periodically forced damped pendulum \eqref{eqn:fdp} in Section~4.
 
% **********************************************************
% Section 2
% **********************************************************

\section{Real-Analytic Nonintegrability}

In this section we review the theory of \cite{MY21b}
 for forced nonlinear oscillators of the form \eqref{eqn:syse}.
See \cite{MY21b} for more details including proofs of the theorems.

\begin{figure}[t]
\begin{center}
\includegraphics[scale=1]{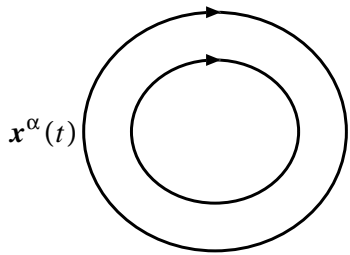}
\caption{Assumption~(M1).\label{fig:2a}}
\end{center}
\end{figure}

We first make the following assumptions on the unperturbed system \eqref{eqn:sys0}:
\begin{enumerate}
\setlength{\leftskip}{-0.8em}
\item[\bf(M1)]
There exists a one-parameter family of periodic orbits $x^\alpha(t)$, $\alpha\in(\alpha_1,\alpha_2)$,
 with period $T^\alpha>0$ for some $\alpha_1<\alpha_2$ (see Fig.~\ref{fig:2a});
\item[\bf(M2)]
$x^\alpha(t)$ is analytic with respect to $\alpha\in(\alpha_1,\alpha_2)$.
\end{enumerate}
Note that in (M1)
 $x^\alpha(t)$ is automatically analytic with respect to $t$
 since the vector field of \eqref{eqn:sys0} is analytic.
Letting $\theta=\omega t \mod 2\pi$ such that $\theta\in\Sset^1$,
 we rewrite \eqref{eqn:syse} as an autonomous system, 
\begin{equation}
\dot{x}=J\D H(x)+\epsilon g(x,\omega t),\quad
\dot{\theta}=\omega.
\label{eqn:asyse}
\end{equation}

We assume that at $\alpha=\alpha^{m/n}$
\begin{equation}
\frac{2\pi}{T^\alpha}=\frac{n}{m}\omega,
\label{eqn:res4}
\end{equation}
where $m,n>0$ are relatively prime integers.
We define the \emph{subharmonic Melnikov function} as
\begin{equation}
M^{m/n}(\theta)=\int_0^{2\pi m/\omega}\D H(x^{\alpha}(t))\cdot g(x^\alpha(t),\omega t+\theta)\d t,
\label{eqn:subM}
\end{equation}
where $\alpha=\alpha^{m/n}$.
If $M^{m/n}(\theta)$ has a simple zero at $\theta=\theta_0$
 and $\d T^\alpha/\d\alpha\neq 0$ at $\alpha=\alpha^{m/n}$,
 then for $|\epsilon|>0$ sufficiently small
 there exists a $2\pi m/\omega$-periodic orbit
 near $(x,\theta)=(x^\alpha(t),\omega t+\theta_0)$ in \eqref{eqn:asyse}.
See Theorem~3.1 of \cite{Y96}.
A similar result is also found in \cite{GH83a,GH83,W90}.
The stability of the periodic orbit can also be determined easily \cite{Y96}.
Moreover, several bifurcations of periodic orbits
 when $\d T^\alpha/\d\alpha\neq 0$ or not
 were discussed in \cite{Y96,Y02,Y03}.

We now state the results of \cite{MY21b}.
Choose a point $\alpha=\alpha_0\in(\alpha_1,\alpha_2)$
 such that $\d T^\alpha/\d\alpha\neq 0$,
 and let $U\subset(\alpha_1,\alpha_2)$ be a neighborhood of $\alpha_0$.
Let
\[
D_\R=\{\alpha^{m/n}\mid \text{$m,n\in\Nset$ are relatively prime}\}\cap U.
\]
A subset $D\subset U$ is called a \emph{key set} for $C^\omega(U)$
 if any analytic function vanishing on $D$ vanishes on $U$.
In particular, if $D$ has an accumulation point in $U$,
 then it becomes a key set for $C^\omega(U)$.
We have the following two theorems.

\begin{thm}\label{thm:2a}
Suppose that there exists a key set $D\subset D_\R$ for $C^\omega(U)$ 
 such that $M^{m/n}(\theta)$ is not identically zero for $\alpha^{m/n}\in D$.
Then %for $|\epsilon|\neq 0$ sufficiently small
 the system \eqref{eqn:asyse} has no real-analytic first integral
 in a neighborhood of $\{x^{\alpha_0}(t)\mid t\in[0,T^{\alpha_0})\}\times\Sset^1$
 such that it depends real-analytically on $\epsilon$ near $\epsilon=0$.
\end{thm}

\begin{thm}\label{thm:2b}
Suppose that there exists a key set $D\subset D_\R$ for $C^\omega(U)$ 
 such that $M^{m/n}(\theta)$ is not constant for $\alpha^{m/n}\in D$.
Then %for $|\epsilon|\neq 0$ sufficiently small
 the system \eqref{eqn:asyse} is not real-analytically Bogoyavlenskij-integrable
 in a neighborhood of $\{x^{\alpha_0}(t)\mid t\in[0,T^{\alpha_0})\}\times\Sset^1$
 such that the first integrals and commutative vector fields also depend analytically
 on $\epsilon$ near $\epsilon=0$.
\end{thm}

\begin{figure}
\begin{center}
\includegraphics[scale=0.9]{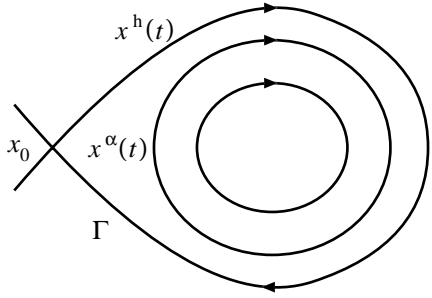}
\caption{Assumption~(M3).\label{fig:M3}}
\end{center}
\end{figure}

We additionally assume the following on the unperturbed system \eqref{eqn:sys0}:
\begin {enumerate}
\setlength{\leftskip}{-0.6em}
\item[\bf(M3)]
There exists a hyperbolic saddle $x_0$ with a homoclinic orbit $x^\h(t)$ such that
\[
\lim_{\alpha\to\alpha_2}\sup_{t\in\Rset}d(x^\alpha(t),\Gamma)=0,
\]
where $\Gamma=\{x^\h(t)\mid t\in\Rset\}\cup\{x_0\}$
 and $d(x,\Gamma)=\inf_{y\in\Gamma}|x-y|$.
See Fig.~\ref{fig:M3}.
\end{enumerate}
We define the \emph{homoclinic Melnikov function} as 
\begin{equation}
M(\theta)=\int_{-\infty}^{\infty}\D H(x^{\h}(t))\cdot g(x^{\h}(t),t+\theta)\d t.
\label{eqn:homM}
\end{equation}
If $M(\theta)$ has a simple zero,
 then for $|\epsilon|>0$ sufficiently small
 there exist transverse homoclinic orbits to a periodic orbit near $\{x_0\}\times\Sset^1$
 in \eqref{eqn:asyse} \cite{GH83,M63,W90}.
The existence of such transverse homoclinic orbits
 implies that the system \eqref{eqn:asyse} exhibits chaotic motions
 by the Smale-Birkhoff theorem \cite{GH83,M73,W90}
 and has no (additional) real-analytic first integral (see, e.g., Theorem~3.10 of \cite{M73}).
We easily show that
\begin{equation}
\lim_{m\to \infty} M^{m/1}(\theta)=M(\theta)
\label{eqn:limM}
\end{equation}
for each $\theta\in \Sset^1$ (see Theorem~4.6.4 of \cite{GH83}).
Let $U$ be a neighborhood of $\alpha=\alpha_2$.
It follows from \eqref{eqn:limM} that
 if $M(\theta)$ is not identically zero or constant,
 then for $m>0$ sufficiently large neither is $M^{m/1}(\theta)$.
Let $\hat{U}\subset\Rset^2$ be a region such that $\partial\hat{U}\supset \Gamma$
%where $\Gamma^\p=\{x^\alpha(t)\mid t\in[0,\hat{T}^\alpha] \}$
and $\hat{U}\supset\{x^\alpha(t)\mid t\in[0,T^\alpha)\}$
for some $\alpha\in(\alpha_1,\alpha_2)$.
Noting the relation \eqref{eqn:limM},
 we obtain the following from Theorems~\ref{thm:2a} and \ref{thm:2b}.

\begin{thm}\label{thm:2c}
Suppose that $M(\theta)$ is not identically zero.
Then the system \eqref{eqn:asyse} has no real-analytic first integral in $\hat{U}\times\Sset^1$
 such that it depends real-analytically on $\epsilon$ near $\epsilon=0$.
\end{thm}

\begin{thm}\label{thm:2d}
Suppose that $M(\theta)$ is not constant.
Then the system \eqref{eqn:asyse} is not real-analytically Bogoyavlenskij-integrable
 in $\hat{U}\times\Sset^1$
 such that the first integrals and commutative vector fields also depend analytically
 on $\epsilon$ near $\epsilon=0$.
\end{thm}

\begin{rmk}
\label{rmk:2a}
Theorems~$\ref{thm:2c}$ and $\ref{thm:2d}$, respectively,
 mean that the system \eqref{eqn:asyse} may have no first integral and be Bogoyavlenskij-nonintegrable
 even if the Melnikov function $M(\theta)$ does not have a simple zero,
 i.e., there may exist no transverse homoclinic orbit to the periodic orbit,
 but it is not identically zero and constant.
See also Section~$6.3$ of {\rm\cite{MY21b}} and Remark~$\ref{rmk:4a}$.
\end{rmk}

% **********************************************************
% Section 3
% **********************************************************

\section{Complex-Meromorphic Nonintegrability}

In this section we review the theory of \cite{Y21b}
 for forced nonlinear oscillators of the form \eqref{eqn:syse}.
See \cite{Y21b} for more details including a proof of the theorem.

\begin{figure}
\begin{center}
\includegraphics[scale=0.8]{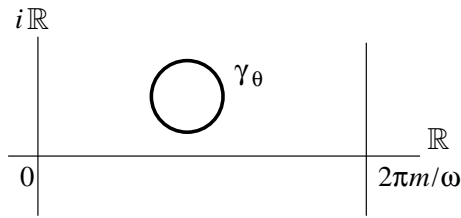}
\caption{Closed path $\gamma_\theta$.\label{fig:3a}}
\end{center}
\end{figure}

We assume (M1) and (M2) on the unperturbed system \eqref{eqn:sys0}.
For \eqref{eqn:asyse}
 we extend the domain of the independent variable $t$ to a domain including $\Rset$ in $\Cset$
 and do so for the dependent variables $x$ and $\theta$.
Let $\gamma_\theta$ be a closed path in a domain containing $(0,2\pi m/\omega)$ in $\Cset$
 such that $\gamma_\theta\cap(i\Rset\cup(2\pi m/\omega+i\Rset))=\emptyset$.
See Fig.~\ref{fig:3a}.
We assume that at $\alpha=\alpha^{m/n}$ the resonance condition~\eqref{eqn:res4} holds
 for $m,n>0$ relatively prime integers, as in Section~2.
For $\alpha=\alpha^{m/n}$ we define the integral
\begin{equation}
\hat{\I}(\theta)
=\int_{\gamma_\theta}\D H(x^\alpha(t))
 \cdot g\left(x^\alpha(t),\omega t+\theta\right)\d t,
\label{eqn:thm3a}
\end{equation}
which is similar to the subharmonic Melnikov function \eqref{eqn:subM}
 but defined by a contour integral along the closed circle $\gamma_\theta$.
An integral similar to \eqref{eqn:thm3a} for not periodic but homoclinic orbits
 was used in \cite{M02,Z82a}.
We have the following theorem.

\begin{thm}
\label{thm:3a}
Suppose that at $\alpha=\alpha^{m/n}$, $\d T^\alpha/\d\alpha\neq 0$
 and there exists a closed loop $\gamma_\theta$
 in a domain including $(0,2\pi m/\omega)$ in $\Cset$
 such that $\gamma_\theta\cap(i\Rset\cup(2\pi m/\omega+i\Rset))=\emptyset$
 and the integral $\hat{\I}(\theta)$ is not zero for some $\theta=\theta_0\in\Sset^1$.
Then the system \eqref{eqn:asyse}
 is not complex-meromorphically Bogoyavlenskij-integrable
 near the resonant periodic orbit $(x,\theta)=(x^\alpha(t),\omega t+\theta_0)$
 with $\alpha=\alpha^{m/n}$ on any domain $\hat{\Gamma}$ in $\Cset/(2\pi m/\omega)\Zset$
 containing $\Rset/(2\pi m/\omega)\Zset$ and $\gamma_\theta$,
 such that the first integrals and commutative vector fields also depend complex-meromorphically
 on $\epsilon$ near $\epsilon=0$.
Moreover, if the integral $\hat{\I}(\theta)$ is not zero for any $\theta\in\hat{\Delta}$,
 where $\hat{\Delta}$ is a dense set of $\Sset^1$,
 then the conclusion holds for any periodic orbit on the resonant torus
 $\T^\ast=\{(x^\alpha(t),\omega t+\theta)\mid t\in\hat{\Gamma},\theta\in\Sset^1,\alpha=\alpha^{m/n}\}$.
\end{thm}

\begin{rmk}
\label{rmk:3a}
Theorem~$\ref{thm:3a}$
 means that the system \eqref{eqn:asyse} may be Bogoyavlenskij-nonintegrable
 even if the Melnikov function $M(\theta)$ does not have a simple zero,
 i.e., there may exist no transverse homoclinic orbit to the periodic orbit.
See also Section~$4.3$ of {\rm\cite{Y21b}} and Remark~$\ref{rmk:4a}$.
\end{rmk}

% **********************************************************
% Section 4
% **********************************************************

\section{Forced Damped Pendulum}

\begin{figure}
\begin{center}
\includegraphics[scale=0.25]{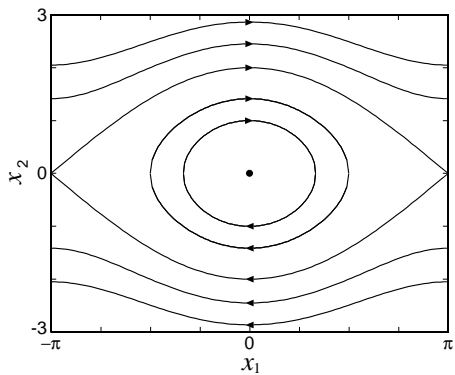}
\caption{Phase portraits of \eqref{eqn:fdp} with $\epsilon=0$.\label{fig:4a}}
\end{center}
\end{figure}

We apply the above two theories
 and discuss the Bogoyavlenskij-nonintegrability
 of the periodically forced damped pendulum \eqref{eqn:fdp}.
 
Letting $\theta=\omega t \mod 2\pi$, we rewrite  \eqref{eqn:fdp} as an autonomous system
\begin{equation}
\dot{x}_1=x_2,\quad
\dot{x}_2=-\sin x_1+\epsilon(\beta\cos\omega t-\delta x_2),\quad
\dot{\theta}=\omega,
\label{eqn:afdp}
\end{equation}
as in \eqref{eqn:asyse}.
When $\epsilon=0$,
 Eq.~\eqref{eqn:fdp} becomes a single-degree-of-freedom Hamiltonian system \eqref{eqn:sys0}
 with the Hamiltonian
\[
H=1-\cos x_1+\frac{1}{2}x_2^2,
\]
and has an elliptic point at $(0,0)$ and a hyperbolic saddle at $(\pi,0)$.
The phase portraits of \eqref{eqn:fdp} with $\epsilon=0$ are shown in Fig.~\ref{fig:4a}:
 There exist a pair of homoclinic orbits
\[
x^\h_\pm(t)=\bigl(\pm2\arcsin(\tanh t), \pm2\sech t\bigr),
\]
a one-parameter family of periodic orbits
\begin{align*}
x^k(t)
 =&\bigl(2\arcsin(k\sn t), 2k\cn t\bigr),\quad
k\in(0,1),
\end{align*}
between the homoclinic orbits,
 and a pair of one-parameter families of periodic orbits
\begin{align*}
\tilde{x}^k_{\pm}(t)
 =&\biggl(\pm 2\arcsin\left(\sn \left(\frac{t}{k}\right)\right), \pm \frac{2}{k} \dn \left(\frac{t}{k}\right)\biggr),\quad
k\in(0,1),
\end{align*}
above or below the homoclinic orbits,
where $\sn$, $\cn$ and $\dn$ represent the Jacobi elliptic functions with the elliptic modulus $k$.
The periods of $x^k(t)$ and $\tilde{x}^k_{\pm}(t)$ are given
 by $\hat{T}^k=4K(k)$ and $\tilde{T}^k=2kK(k)$, respectively,
 where $K(k)$ is the complete elliptic integral of the first kind.
Note that $\tilde{x}^k_{\pm}(t)$ approaches $x^\h_{\pm}(t)$ as $k\to 1$.
The above analytical expressions of these orbits are found, e.g., in \cite{BCMSP87,BW99}.
See \cite{BF54,WW27} for general information on elliptic functions.

We first apply the theory of Section~2.
Assume that the resonance conditions
\begin{equation}
n\hat{T}^k=\frac{2\pi m}{\omega},\quad\mbox{i.e.,}\quad
\omega=\frac{\pi m}{2nK(k)},
\label{eqn:resk}
\end{equation}
and
\begin{equation}
n\tilde{T}^k=\frac{2\pi m}{\omega},\quad\mbox{i.e.,}\quad
\omega=\frac{\pi m}{nkK(k)},
\label{eqn:tresk}
\end{equation}
hold for $x^k(t)$ and $\tilde{x}^k_\pm(t)$, respectively,
 with $m,n>0$ relatively prime integers.
We compute the subharmonic Melnikov function \eqref{eqn:subM}
 for $x^k(t)$ and $\tilde{x}^k_\pm(t)$ as
\begin{align}\label{eqn:M1}
M^{m/n}(\theta)=
-\delta J_1(k,m)+\beta J_2(k,m,n)\cos\theta
%\begin{cases}
%-16n(E(k)-k'^2K(k))\delta+\frac{2\pi}{k}\sech\left(\frac{m\pi K(k')}{2K(k)}\right) & \mbox{(for $n=1$ and $m$ odd)};\\
%-16n(E(k)-k'^2K(k))\delta \quad & \mbox{(for $n \neq 1$ or $m$ even)},\
%\end{cases}
\end{align}
and
\begin{align}\label{eqn:M2}
\tilde{M}^{m/n}_\pm(\theta)=-\delta\tilde{J}_1(k,m)\pm\beta\tilde{J}_2(k,m,n)\cos\theta,
\end{align}
respectively, where
\begin{align*}
&
J_1(k,n)=16n\bigl( E(k)-k'^2K(k)\bigr),\\
&
J_2(k,m,n)=
\begin{cases}
%4\pi\sech\left(\displaystyle\frac{m\pi K(k')}{2K(k)}\right)
4\pi\sech\bigl(\omega K(k')\bigr)
 & \mbox{(for $n=1$ and $m$ odd)}; \\
0 & \mbox{(for $n\neq 1$ or $m$ even).}
\end{cases}
\\
&
\tilde{J}_1(k,n)=\frac{8nE(k)}{k},\\
&
\tilde{J}_2(k,m,n)=\begin{cases}
%2\pi \sech\left(\displaystyle\frac{m\pi K(k')}{K(k)}\right)
2\pi \sech\bigl(k\omega K(k')\bigr)
 & \mbox{(for $n=1$)};\\
0\quad & \mbox{(for $n\neq 1$)}.
\end{cases}
\end{align*}
Here $E(k)$ is the complete elliptic integral of the second kind 
 and $k'=\sqrt{1-k^2}$ is the complimentary elliptic modulus.
When $\delta>0$,
 the subharmonic Melnikov functions $M^{m/n}(\theta)$ and $\tilde{M}^{m/n}_\pm(\theta)$
 are not identically zero for any relatively prime integers $m,n>0$
 since $J_1(k,n)$ and $\tilde{J}_1(k,n)$ are not zero.
We also compute the homoclinic Melnikov function \eqref{eqn:homM} for $x_\pm^\h(t)$ as
\begin{align}\label{eqn:M3}
M_\pm(\theta)
 =-8\delta\pm2\pi\beta\sech\left(\frac{\pi\omega}{2}\right)\cos\theta.
\end{align}
See \cite{HH86} for more details on the derivation of \eqref{eqn:M2} and \eqref{eqn:M3}. 
Equation~\eqref{eqn:M1} is obtained similarly by the method of residues. 
Note that $M_\pm(\theta)$ has a simple zero when condition \eqref{eqn:chaos} holds.

Let
\begin{align*}
&
R=\{k\in(0,1)\mid \mbox{$k$ satisfies \eqref{eqn:resk} for some $m, n\in\Nset$}\},\\
&
\tilde{R}=\{k\in(0,1)\mid\mbox{$k$ satisfies \eqref{eqn:tresk} for some $m,n\in\Nset$}\},
\end{align*}
and let
\begin{align*}
&
\Gamma_\pm=\{x_\pm^\h(t)\in\Rset^2\mid t\in\Rset\}\cup\{(\pi,0)\},\\
&
S^k=\{(x^k(t),\theta)\in\Rset^2\times\Sset^1\mid t\in[0,\hat{T}^k),\theta\in\Sset^1\},\\
&
\tilde{S}^k_\pm
 =\{(\tilde{x}^k_\pm(t),\theta)\in\Rset^2\times\Sset^1\mid t\in[0,\tilde{T}^k),\theta\in\Sset^1\},\\
&
S^\h_\pm=\Gamma_\pm\times \Sset^1.
\end{align*}
Using Theorems~\ref{thm:2a} and \ref{thm:2c}, we have the following.

\begin{prop}
\label{prop:4a}
The system~\eqref{eqn:afdp} has no real-analytic first integral
depending real-analytically on $\epsilon$ near $\epsilon=0$ in neighborhoods of $S^k$ for $k\in R$,
 of $\tilde{S}^k_\pm$ for $k\in\tilde{R}$, and of $S^h_\pm$ if $\delta>0$.
\end{prop}

Noting that
\[
\lim_{m\to\infty}M^{m/1}_\pm(\theta)=M_+(\theta)+M_-(\theta)
\]
and using Theorem~\ref{thm:2d} and its slight extension, we also have the following.

\begin{prop}
\label{prop:4b}
Let $\hat{U}$ $($resp. $\tilde{U}_\pm)$ be a region $($resp. regions$)$ in $\Rset^2$
 such that $\partial\hat{U}\supset\Gamma_\pm$ $($resp. $\partial\tilde{U}_\pm\supset\Gamma_\pm)$
 and $\hat{U}\supset\{x^k(t)\mid t\in[0,\hat{T}^k) \}$
 $($resp. $\tilde{U}_\pm\supset\{\tilde{x}_\pm^k(t)\mid t\in[0,\tilde{T}^k) \})$
for some $k\in(0,1)$.
If $\beta>0$, then the system~\eqref{eqn:afdp} is not real-analytically Bogoyavlenskij-integrable
 in the meaning of Theorem~$\ref{thm:2d}$
 in $\hat{U}\times\Sset^1 $ $($resp. in $\tilde{U}_\pm\times\Sset^1)$.
\end{prop}

We next apply the theory of Section~3.
The integral \eqref{eqn:thm3a} becomes
\begin{equation}
\hat{\I}^k(\theta)
 =-4k^2\delta\int_{\gamma_\theta}\cn^2t\,\d t
 +2k\beta\int_{\gamma_\theta}\cn t\,\cos(\omega t+\theta)\d t
\label{eqn:int10}
\end{equation}
and
\begin{equation}
\hat{\I}_\pm^k(\theta)
 =-\frac{4\delta}{k^2}\int_{\gamma_\theta}\dn^2\left(\frac{t}{k}\right)\d t
 +\frac{2\beta}{k}\int_{\gamma_\theta}\dn\left(\frac{t}{k}\right)\cos(\omega t+\theta)\d t
\label{eqn:int20}
\end{equation}
for $x^k(t)$ and $\tilde{x}_\pm^k(t)$, respectively.
We shift the variable $t$ by $\frac{1}{2}\hat{T}^k$ (resp. $\frac{1}{2}\tilde{T}^k$)
 in \eqref{eqn:int10} (resp. in \eqref{eqn:int20}),
 and take a circle centered at $t=iK(k')+\frac{1}{2}\hat{T}^k$
 (resp. $t=ikK(k')+\frac{1}{2}\tilde{T}^k$) with sufficiently small radius,
 as $\gamma_\theta$ (resp. $\tilde{\gamma}_\theta$).
So we compute \eqref{eqn:int10} and \eqref{eqn:int20}, respectively, as 
\begin{align}
\hat{\I}^k(\theta)=4\pi\beta\biggl(\cosh \bigl(\omega K(k')\bigr)\cos\theta-i\sinh \bigl(\omega K(k')\bigr)\sin \theta\biggr)
\label{eqn:int1}
\end{align}
and
\begin{align}
\hat{\I}^k_\pm(\theta)
 =\pm4\pi\beta\biggl(\cosh \bigl(\omega k K(k')\bigr)\cos\theta
 -i\sinh \bigl(\omega k K(k')\bigr)\sin \theta\biggr),
\label{eqn:int2}
\end{align}
which are not zero for any $\theta\in\Sset^1$ if $\beta>0$.
See Appendix~\ref{app:deriv-cpx} for the derivation of \eqref{eqn:int1} and \eqref{eqn:int2}.

Let $\hat{\Gamma}$ (resp. $\tilde{\Gamma}$) be a domain in $\Cset/(2\pi m/\omega)\Zset$
 containing $\Rset/(2\pi m/\omega)\Zset$ and $t=iK(k')+\frac{1}{2}\hat{T}^k$
 (resp. $t=ikK(k')+\frac{1}{2}\tilde{T}^k$).
For $k\in R$ and $k\in \tilde{R}$, let
\[
\T^k=\{(x^k(t),\omega t+\theta)\mid t\in \hat{\Gamma},\theta\in\Sset^1\},
\]
and
\[
\tilde{\T}^k_\pm=\{(\tilde{x}^k_\pm(t),\omega t+\theta)\mid t\in \tilde{\Gamma},\theta\in\Sset^1\},
\]
respectively.
Using Theorem~\ref{thm:3a}, we have the following.
\begin{prop}\label{prop:4c}
If $\beta>0$, then the system \eqref{eqn:afdp}
 is complex-meromorphically Bogoyavlenskij-nonintegrable
 in the meaning of Theorem~$\ref{thm:3a}$ near any periodic orbit
 on $\T^k$ with $k\in R$ and on $\tilde{\T}^k_\pm$ with $k\in \tilde{R}$.
\end{prop}

\begin{rmk}\
\label{rmk:4a}
\begin{enumerate}
\setlength{\leftskip}{-1.8em}
\item[\rm(i)]
Proposition~$\ref{prop:4a}$ shows that
 when $\delta>0$, the system \eqref{eqn:fdp}
 has no first integral depending real-analytically on $\epsilon$ near $\epsilon=0$
 even if condition~\eqref{eqn:chaos} does not hold,
 i.e., there exists no transverse homoclinic orbit near $\Gamma_\pm$.
\item[\rm(ii)]
Propositions~$\ref{prop:4b}$ and $\ref{prop:4c}$, respectively, show that
 when $\beta>0$, the system \eqref{eqn:fdp}
 is Bogoyavlenskij-nonintegrable in the meanings of Theorems~$\ref{thm:2d}$ and $\ref{thm:3a}$,
 even if condition~\eqref{eqn:chaos} does not hold.
\item[\rm(iii)]
If $\beta=0$, then Propositions~$\ref{prop:4b}$ and $\ref{prop:4c}$ say nothing
 about the nonintegrability of \eqref{eqn:afdp}.
However, if $\delta>0$ holds additionally, then by Proposition~$\ref{prop:4a}$
 the system \eqref{eqn:afdp} has no real-analytic first integral
 depending on $\epsilon$ analytically near $\epsilon=0$.
\end{enumerate}
\end{rmk}

% **********************************************************
% Appendices
% **********************************************************

\appendix
\renewcommand{\theequation}{\Alph{section}.\arabic{equation}}
%\setcounter{equation}{0}

%\section{Derivation of \eqref{eqn:M1}, \eqref{eqn:M2} and \eqref{eqn:M3}} \label{app:deriv-real}

\section{Derivation of \eqref{eqn:int1} and \eqref{eqn:int2}}\label{app:deriv-cpx}
In this appendix we use the method of residues
 and compute the integrals \eqref{eqn:int10} and \eqref{eqn:int20}.
A similar calculation is found in \cite{Y21b}.

%For the system \eqref{eqn:afdp}, the integral \eqref{eqn:thm3a} becomes
%\begin{align}
%\label{eqn:int-fdp}
%	\hat{\I}(\theta)=-\delta\int_{\gamma_\theta} (x^\alpha(t))^2 \d t
%				&+\beta\left(\int_{\gamma_\theta} x^\alpha(t)\cos(\omega t) \d t\right)\cos\theta\nonumber\\
%				&\qquad-\beta\left(\int_{\gamma_\theta} x^\alpha(t)\sin(\omega t) \d t\right)\sin \theta.
%\end{align}
We begin with the first term in \eqref{eqn:int10}.
Letting $s=1/\sn t$, we have
\begin{align}\label{eqn:int-a}
\int \cn^2 t\, \d t=-\int \frac{1}{s^2} \sqrt{\frac{1-s^2}{k^2-s^2}}\d s
\end{align}
from the basic properties of the Jacobi elliptic functions
\[
\frac{\d}{\d t} \sn t=\cn t \dn t,
\quad \cn^2 t=1-\sn^2 t,
\quad \dn^2 t=1-k^2\sn^2 t.
\]
Obviously, the integrand in the right hand side of \eqref{eqn:int-a} has a pole of order $2$ and
\[
\frac{\d}{\d s}\sqrt{\frac{1-s^2}{k^2-s^2}}=0
\]
at $s=0$.
Noting that $s=0$ when $t=iK(k')$, we obtain
\begin{align}\label{eqn:int-b}
\int_{\bar{\gamma}_\theta}\cn^2 t\, \d t
=\int_{|s|=\rho} \frac{1}{s^2} \sqrt{\frac{1-s^2}{k^2-s^2}}\d s=0
\end{align}
by the method of residues,
 where $\bar{\gamma}_\theta=\{t\in\Cset\mid t+\tfrac{1}{2}\hat{T}^k\in\gamma_\theta\}$
 and $\rho>0$ is sufficiently small.

We turn to the second term in \eqref{eqn:int10}.
%For constants $A,B\in \Cset$, we compute
%\begin{align}\label{eqn:cos}
%\cos At&=\cos(A(t-iB)+iAB)\nonumber\\
%&=\cosh (AB) \cos A(t-iB)-i\sinh(AB)\sin A(t-iB)
%\end{align}
%from
%\[
%\cos(it)=\cosh t,\quad \sin(it)=i\sinh t.
%\]
%Eq. \eqref{eqn:cos} implies that
%near $t=iK(k')$.
Since
\begin{align}\label{eqn:cos-ex}
\cos \omega t=\cosh (\omega K(k')) +O(t-iK(k'))
\end{align}
and
\[
\cn t=-\frac{i}{k(t-iK(k'))}+O(1)
\]
near $t=iK(k')$, we have
%\[
%\cn t\, \cos \omega t=-\frac{i \cosh\omega K(k')}{k(t-iK(k'))}+O(1)
%\]
%there, so that
\[
\int_{\bar{\gamma}_\theta} \cn t\cos \omega t\,\d t=\frac{2\pi}{k} \cosh \omega K(k').
\]
Similarly, since
\begin{align}
\sin \omega t
%&=i\sinh (\omega K(k')) \cos \omega (t-iK(k'))+\cosh(\omega K(k'))\sin \omega (t-iK(k'))\nonumber\\
&=i\sinh (\omega K(k'))+O(t-iK(k'))\label{eqn:sin}
\end{align}
near $t=iK(k')$, we have
\[
\int_{\bar{\gamma}_\theta} \cn t\sin \omega t\,\d t=\frac{2\pi i}{k} \sinh\omega K(k').
\]
Thus, we obtain \eqref{eqn:int1}.

We next compute \eqref{eqn:int20}.
We easily see that the first term vanishes by \eqref{eqn:int-b} since
\[
\dn^2 t%=1-k^2\sn^2 \zeta=1-k^2(1-\cn^2 t)
=k'^2-k^2\cn^2 t.
\]
On the other hand, since
\begin{align*}
\dn t=-\frac{i}{t-iK(k')}+O(1),
\end{align*}
we have
\begin{align*}
\int_{\bar{\gamma}_\theta} \dn\left(\frac{t}{k}\right)\cos \omega t\,\d t
=2\pi k \cosh(\omega kK(k'))
\end{align*}
by \eqref{eqn:cos-ex}, 
 where $\bar{\gamma}_\theta=\{t\in\Cset\mid t+\tfrac{1}{2}\tilde{T}^k\in\gamma_\theta\}$.
 Similarly, by \eqref{eqn:sin} we have
\begin{align*}
\int_{\bar{\gamma}_\theta} \dn\left(\frac{t}{k}\right)\sin \omega t\,\d t
=2\pi ik \sinh (\omega kK(k')).
\end{align*}
Thus, we obtain \eqref{eqn:int2}.

% **********************************************************
% Bibliography
% **********************************************************

\end{document}